\documentclass[preprint,1p,10pt]{elsarticle}

\journal{arXiv}

\biboptions{sort&compress,comma,round}

\usepackage{graphicx}
\usepackage{amssymb}
\usepackage[FIGTOPCAP]{subfigure}
\usepackage{stmaryrd}
\usepackage{algorithm}
\usepackage{algpseudocode} 

\usepackage{color}

\usepackage{graphics}
\usepackage{amssymb, latexsym, amsmath, epsfig}
\usepackage{graphicx}
 \usepackage{color}
 \usepackage{epstopdf}

\numberwithin{equation}{section}

\newtheorem{remark}{Remark}

\newcommand{\bfs}[1]{{\boldsymbol #1}}

\begin{document}

\begin{frontmatter}


\title{Split generalized-$\alpha$ method: A linear-cost solver for a modified generalized-method for multi-dimensional second-order hyperbolic systems}

\author[ad1]{Pouria Behnoudfar\corref{corr}}
\ead{pouria.behnoudfar@postgrad.curtin.edu.au}

\author[ad1,ad2]{Quanling Deng}
\cortext[corr]{Corresponding author}
\ead{Quanling.Deng@curtin.edu.au}
\author[ad1,ad2]{Victor M. Calo}
\ead{Victor.Calo@curtin.edu.au}


\address[ad1]{School of Earth and Planetary Sciences, Curtin University, Kent Street, Bentley, Perth, WA 6102, Australia}

\address[ad2]{Mineral Resources, Commonwealth Scientific and Industrial Research Organisation (CSIRO), Kensington, Perth, WA 6152, Australia}

\begin{abstract}
We propose a variational splitting technique for the generalized-$\alpha$ method to solve hyperbolic partial differential equations. We use tensor-product meshes to develop the splitting method, which has a computational cost that grows linearly with respect to the total number of degrees of freedom for multi-dimensional problems. We consider standard $C^0$ finite elements as well as smoother B-splines in isogeometric analysis for the spatial discretization. We also study the spectrum of the amplification matrix to establish the unconditional stability of the method. We then show that the stability behavior affects the overall behavior of the integrator on the entire interval and not only at the limits $0$ and $\infty$. We use various examples to demonstrate the performance of the method and the optimal approximation accuracy. For the numerical tests, we compute the $L_2$ and $H^1$ norms to show the optimal convergence of the discrete method in space and second-order accuracy in time. 
\end{abstract}


\end{frontmatter}

\vspace{0.5cm}

\section{Introduction}

The generalized-$\alpha$ method is an implicit method for solving dynamic problems, which is second-order accurate and provides user-controlled numerical dissipation. Chung and Hulbert in~\cite{ chung1993time} proposed this method to solve structural dynamics problems in which the underlying partial differential equations are hyperbolic. The generalized-$\alpha$ method balances the high and low-frequency dissipation. This algorithm minimizes the low-frequency dissipation while it maximizes the high-frequency dissipation~\cite{ chung1993time}.

Various methods provide control over dissipation. For example, the Newmark methods~\cite{ newmark1959method} control dissipation, whereas it has high dissipation in the low-frequency region. Similarly, the $\theta$ method~\cite{ wilson1968computer}, the $\theta_1$ method~\cite{ hoff1988development}, and the $\rho$ method proposed by Bazzi and Anderheggen~\cite{ bazzi1982rho} have high dissipation in the low-frequency region. The generalized-$\alpha$ method includes the description of the aforementioned these time integrators in a single framework. Furthermore, particular choices of the free parameters, reduce the algebraic system resulting from the generalized-$\alpha$ method to those of the HHT-$\alpha$~\cite{ hilber1977improved} and WBZ-$\alpha$~\cite{ wood1980alpha} methods.

Splitting techniques approximate the linear systems and to reduce the computation cost~\cite{ SPORTISSE2000140}. Reducing the dimension of the matrices resulting from implicit numerical schemes, ~\cite{ SPORTISSE2000140} and domain decomposition methods~\cite[Chapter 8]{ SMV02} are two examples.  In this work, we use the idea of operator splitting and adapt it to tensor-product approximations in the context of the classical finite elements and isogeometric analysis for the spatial discretization.  We use splitting methods to reduce the computational cost of this class of meshes significantly. We solve the multi-dimensional problems with computational time and storage that grow proportionally to the number of degrees of freedom in the system. The results reported in~\cite{ gao2014fast, gao2015preconditioners, los2015dynamics, los2017application} show that alternating direction splitting solvers based on tensor-product provide an overall linear computational cost at every time step for various problems. In~\cite{ behnoudfar2018variationally}, we introduced a variational splitting for parabolic problems, and herein, we extend the idea to hyperbolic problems. That is, we propose a splitting for the generalized-$\alpha$ method for hyperbolic equations. We use tensor-product grids to formulate the variational formulation for multi-dimensional problems. We then write the $d$-dimensional formulation as a product of $d$ formulations in each dimension plus error terms. We refer to these formulations as variationally separable. Based on the variational separability, we present a splitting technique to solve the resulting linear systems with a linear computational cost. With sufficient regularity, the approximate solution converges to the exact solution with optimal rates in space and time while reducing the computational cost significantly.

The rest of this paper is as follows. Section~\ref{sec:2} describes the hyperbolic problem under consideration and introduces the particular spatial discretizations to arrive at the matrix formulations of the problem. Section~\ref{sec:3} presents a temporal discretization using the generalized-$\alpha$ method. Therein, we also introduce various splitting methods. Section~\ref{sec:4} establishes the stability of the splitting schemes. We show numerically in Section~\ref{sec:5} that the approximate solution converges optimally to the exact solution. We also verify that the computational cost is linear. Concluding remarks are given in Section~\ref{sec:6}.

\section{Problem statement}\label{sec:2}
Let $\Omega =(0, 1)^d \subset \mathbb{R}^d, d=2,3,$ be an open and bounded domain. We consider the wave propagation problem
\begin{equation} \label{eq:pde}
\begin{aligned}
u_{tt} - \Delta u & = f, \quad \quad \! \qquad (x, t)  \in \Omega \, \times \, ]0, T[, \\
u & = u_D, \qquad \quad (x, t) \in \partial \Omega \, \times \, ]0,T[, \\
u & = d_0, \qquad\quad \, \, (x, t) \in \Omega, \,t=0, \\
\dot u & = v_0, \qquad \quad  \, \,(x, t) \in \Omega, \,t=0. \\
\end{aligned}
\end{equation}
where $\Delta$ is the Laplacian, $f$ is the forcing function, $d_0$ and $v_0$ are the given initial displacement and velocity, respectively. We write the hyperbolic equation~\eqref{eq:pde} in weak form. Below, we first discretize this problem in the spatial domain using isogeometric analysis. 

\subsection{Spatial semi-discretization}

The spatial discretization results in an approximate solution $u_h(\cdot, t)$ in the finite-dimensional space generated by the isogeometric elements for $t\ge 0$ as a solution of an initial value problem written in a system of ordinary differential equations. The matrix resulting from the wave propagation equation becomes:
\begin{equation} \label{eq:pb}
M \ddot{U} + KU = F,
\end{equation}
where $M$ and $K$ are the mass and stiffness matrices, respectively; see, for example,~\cite{ behnoudfar2018variationally} for the detailed derivations. $F$ is the applied-load vector, which is a given function of time. $U$ is the vector of displacement unknowns, and superposed dots indicate temporal differentiation, that is, $\dot{U}$ and $\ddot{U}$ are the velocity and acceleration vectors, respectively.

\subsection{Fully-discrete time-marching scheme}

 The initial boundary-value problem consists of finding a function $U=U(t)$ which satisfies~\eqref{eq:pb} and the initial conditions 
\begin{equation}
U(0) = U_0, \qquad \dot U(0) = V_0,
\end{equation}
 Considering a time marching $0 = t_0 < t_1 < \cdots < t_N = T$, $U_n, V_n, $ and $A_n$ are given approximations to $U(t_n), \dot U(t_n)$, and $\ddot U(t_n)$, respectively. Using a Taylor expansion, linear expressions for $U_{n+1}$ and $V_{n+1}$ are delivered. $A_{n+1}$ is obtained by considering the equation of~\eqref{eq:pb} at the time-step $n+1$.

\section{The generalized-$\alpha$ based splitting method}\label{sec:3}
We control the numerical dissipation in the higher frequency regions by using the generalized-$\alpha$ method~\cite{chung1993time} which we write as
\begin{subequations} \label{eq:galpha}
\begin{align}
\label{eq:galpha1}
U_{n+1} & = U_n + \tau v_n + \frac{\tau^2}{2} A_n + \tau^2 \beta \llbracket A_n \rrbracket, \\
\label{eq:galpha2}
V_{n+1} & = V_n + \tau A_n + \tau \gamma \llbracket A_n \rrbracket, \\
\label{eq:galpha3}
M A_{n+\alpha_m}  & =- K U_{n+\alpha_f}+F_{n+\alpha_f}, 
\end{align}
\end{subequations}
with initical conditions
\begin{equation}
U_0  = U(0), \quad V_0 = V(0),  \quad A_0 = M^{-1} (F_0 - K U_0), \\
\end{equation}
where $F_{n+\alpha_f} = F(t_{n+\alpha_f})$ and 
\begin{equation} \label{eq:param}
\begin{aligned}
W_{n+\alpha_g} & = W_n + \alpha_g \llbracket W_n \rrbracket,\qquad \llbracket W_n \rrbracket & = W_{n+1} - W_n, \quad W = U, V, A, \quad g=m, f.\\
\end{aligned}
\end{equation}
We use $\alpha_m$ rather than $1-\alpha_m$ as originally used in~\cite{ chung1993time}, which is standard usage~(e.g., see~\cite{ bazilevs2008isogeometric}). The method is second-order accurate in time when
\begin{equation}\label{eq:gam}
\gamma  = \frac{1}{2} + \alpha_m - \alpha_f.
\end{equation}
We selected the parameters as
\begin{equation}\label{eq:stab}
\begin{aligned}
\alpha_m & \ge \alpha_f \ge \frac{1}{2},\\
\end{aligned}
\end{equation}
to achieve unconditional stability.  Finally, the user-control of the high-frequency damping requires the use of the following definitions:
\begin{equation}\label{eq:ro}
\begin{aligned}
\beta  &= \frac{1}{4} (1 + \alpha_m - \alpha_f)^2, \\
\alpha_f & = \frac{1}{1+\rho_\infty},\\
\alpha_m & = \frac{2-\rho_\infty}{1+\rho_\infty}, \\
\end{aligned}
\end{equation}
where $\rho_\infty$ defines the ratio of the amplitude of the highest frequency in the system after the first time step. In order to derive our splitting scheme, we perform the following elimination. From ~\eqref{eq:galpha}, \eqref{eq:gam}, and~\eqref{eq:stab}, we have
\begin{subequations}\label{eq:AVU}
\begin{align}
\label{eq:AVU1}
A_{n+\alpha_m} & = A_n + \alpha_m \llbracket A_n \rrbracket, \\
\label{eq:AVU2}
V_{n+\alpha_f} & = V_n + \alpha_f (V_{n+1} - V_n) \nonumber \\ \nonumber
& = V_n + \alpha_f \big( V_n + \tau A_n + \tau \gamma \llbracket A_n \rrbracket - V_n \big) \\ 
& = V_n + \tau \alpha_f  A_n + \tau \alpha_f \gamma \llbracket A_n \rrbracket, \\ 
\label{eq:AVU3}
\nonumber
U_{n+\alpha_f} & = U_n + \alpha_f (U_{n+1} - U_n) \\\nonumber
& = U_n + \alpha_f \left( U_n + \tau V_n + \frac{\tau^2}{2} A_n + \tau^2 \beta \llbracket A_n \rrbracket - U_n \right) \\
& = U_n + \tau \alpha_f  (V_n + \frac{\tau}{2} A_n ) + \tau^2 \alpha_f \beta \llbracket A_n \rrbracket, 
\end{align}
\end{subequations}
which we plug into~\eqref{eq:galpha3} to yield
\begin{equation} \label{eq:gasp}
\begin{aligned}
\big( \alpha_m M + \tau^2 \alpha_f \beta K \big) \llbracket A_n \rrbracket & = F_{n+\alpha_f} - \left( M A_n + K\left(U_n + \tau \alpha_f V_n + \frac{\tau^2}{2} \alpha_f A_n\right) \right) \\
 &= F_{n+\alpha_f} - \Big( M +  \frac{\tau^2}{2} \alpha_f K \Big) A_n - \Big(  \tau \alpha_f K \Big) V_n - K U_n.
\end{aligned}
\end{equation}

Therefore, we rewrite~\eqref{eq:gasp} as
\begin{equation} \label{eq:wave}
\begin{aligned}
\alpha_m G \llbracket A_n \rrbracket & = F_{n+\alpha_f} - \Big( M + \frac{\tau^2}{2} \alpha_f K \Big) A_n - \tau \alpha_f K V_n - K U_n,
\end{aligned}
\end{equation}
where 
\begin{equation}\label{eq:G}
G = M +  \eta K, \qquad \eta = \frac{\tau^2 \alpha_f \beta} {\alpha_m}.
\end{equation}
Here, after solving~\eqref{eq:wave} and obtaining $A_{n+1}$, we calculate $V_{n+1}$ using the second equation in~\eqref{eq:AVU}. We then calculate $U_{n+1}$ using equation ~\eqref{eq:AVU3} and repeat this procedure at each step as time marches forward. The overall procedure only requires solving one matrix problem at each time step; that is, we invert $G$ only once. This matrix has a tensor-product structure, which allows us to develop splitting schemes. We present and analyze our splitting schemes as follows. 

\subsection{Splitting method}

For a Galerkin discretization of~\eqref{eq:pde} on tensor-product meshes in 2D, we have (for the derivation, see~\cite{ behnoudfar2018variationally})
\begin{equation}
\begin{aligned}
M & = M^x \otimes M^y, \\
K & = K^x \otimes M^y  + M^x \otimes K^y ,
\end{aligned}
\end{equation}
where $M^\xi$ and $K^\xi$ with $\xi = x, y$ are the 1D mass and stiffness matrices, respectively. The operator $\otimes$ is a tensor product. Now, we rewrite
\begin{equation}
G = (M^x +  \eta K^x) \otimes (M^y +  \eta K^y)  - T_\eta,
\end{equation}
where 
\begin{equation}
T_\eta =  \eta^2 K^x \otimes K^y  = \mathcal{O}(\eta^2).
\end{equation}
Thus, $T_\eta$ is of order of $\mathcal{O}(\tau^4)$. Now, to propose our method, we approximate $G$ as
\begin{equation}\label{eq:gtil}
\tilde{G}=(M^x +  \eta K^x) \otimes (M^y +  \eta K^y).
\end{equation}
We split the 3D system using a similar argument.

As we show in section 4, using only $\tilde{G}$ on the left-hand side of~\eqref{eq:wave} lacks the user-control on the numerical dissipation and unconditional stability. Thus we substitute $K$ by an approximation derived from~\eqref{eq:gtil}. To obtain it, we use~\eqref{eq:G} to obtain this equivalent expression:
\begin{equation}
\begin{aligned}
K = \dfrac{1}{\eta} (M+\eta K)-\dfrac{1}{\eta} M.
\end{aligned}
\end{equation}
Then, the approximation we use becomes
\begin{equation}
K \approx \dfrac{1}{\eta} \left(\tilde{G}-M\right).
\end{equation}
Now, we rewrite~\eqref{eq:wave} as
	\begin{equation} \label{spli3}
\begin{aligned}
\alpha_m {\tilde{G}} \llbracket A_n \rrbracket  &= F_{n+\alpha_f} - M A_n-\dfrac{1}{\eta}\left(\tilde{G}-M\right)\cdot \left[\frac{\tau^2\alpha_f }{2}A_n+\tau \alpha_f V_n+U_n\right].
\end{aligned}
\end{equation}

This modification is second-order accurate in time, while it reduces the matrix-assembly cost for the system as only $\tilde{G}$ is needed. Additionally, the system in~\eqref{spli3} provides the user-control on the numerical dissipation. We analyze this method in section~\ref{stb}. Algorithm~\ref{alg:iterative} summarizes the method .
\begin{algorithm}
	\caption{Variationally separable splitting for the generalized-$\alpha$ method for 2D hyperbolic problems}
	\label{alg:iterative}
	\begin{algorithmic}
		\State Set $ T, \tau, \rho_\infty, {M}^\xi, {K}^\xi,{U}_{0}, {V}_{0},
	{F}_{0}$, where $M^\xi$ and $K^\xi$ are the mass an stiffness matrices in $\xi$ direction, where $\xi=x,y$.
		\smallskip 	\smallskip
			\For{$n=1$ until $n=\frac{T}{\tau}$},
				\smallskip 	\smallskip

		\State Set $R_n = F_{n+\alpha_f} - M A_n-\dfrac{1}{\eta}\left(\tilde{G}-M\right)\cdot \left[\frac{\tau^2\alpha_f }{2}A_n+\tau \alpha_f V_n+U_n\right]$.
	
			\smallskip 		\smallskip
		\State Reassemble	$R_n$ as a matrix representation, which rows and columns belong to values corresponding to $x$ and $y$ directions.   
		
				\smallskip 		\smallskip
		\State Solve $\alpha_m (M^x+\eta K^x)\llbracket \tilde{A}_{n} \rrbracket = R_n$ for $\llbracket \tilde{A}_{n} \rrbracket$
		\smallskip 	\smallskip
	\State Solve $(M^y+\eta K^y) \llbracket {A}_{n} \rrbracket =\llbracket \tilde{A}_{n} \rrbracket $ for $\llbracket{A}_{n} \rrbracket$ 
		\smallskip	\smallskip 
		\State Update $V_{n+1}$ from $A_{n+1}$ using~\eqref{eq:AVU2}.
				\smallskip 		\smallskip
				\State Update $U_{n+1}$ form $V_{n+1}$ and $A_{n+1}$ using~\eqref{eq:AVU3}.
					\smallskip 		\smallskip
		\EndFor
			\smallskip
	\end{algorithmic}
\end{algorithm}

\begin{remark}
    The solution cost of the linear system with a matrix $\big(M_\xi+\eta K_\xi\big)^{-1}$, where $\xi=x,y,z$, is linear with respect to the number of degrees of freedom (see for example~\cite{ los2015dynamics, los2017application}). Therefore, the overall cost to solve the resulting multi-dimensional system is linear as it consists of two or three linear systems for 2D and 3D problems, respectively.  
\end{remark}

In the next section, we analyze the method focusing on the stability property and present its advantages. 

\section{Spectral analysis} \label{sec:4}

We start the analysis with the standard generalized-$\alpha$ method. The following analysis follows closely the approach introduced in~\cite{behnoudfar2018variationally, chung1993time}. Throughout this section, we set $F=0$.

\subsection{The generalized-$\alpha$ method}
The standard generalized-$\alpha$ method, combined with~\eqref{eq:wave},~\eqref{eq:galpha1} and~\eqref{eq:galpha2}, results in the following:
\begin{equation} \label{eq:mp0}
\begin{aligned}
\begin{bmatrix}
U^{n+1} \\
\tau V^{n+1}\\
\tau^2 A^{n+1}\\
\end{bmatrix}
& =
\begin{bmatrix}
I - \tau^2 {\alpha_f}^2 J &  I-\tau^2\alpha_f^3 J& \frac{1}{2}I-\alpha_f^2 R \\[0.2cm]
-\tau^2 \gamma J &I-\tau^2\gamma\alpha_fJ &  I-\gamma R \\[0.2cm]
-\tau^2 J&-\tau^2 \alpha_f J&I-R\\
\end{bmatrix}
\cdot\begin{bmatrix}
U^n \\
\tau V^n\\
\tau^2 A^{n}\\
\end{bmatrix},
\end{aligned}
\end{equation}
where $I$ is an identity matrix and we denote
\begin{equation}
J=(\alpha_mG)^{-1} K, \qquad \qquad R=\left(\alpha_mG\right)^{-1}\left(M+\frac{\tau^2}{2}\alpha_fK\right).
\end{equation}
To study its stability, we spectrally decompose  the matrix $K$ with respect to $M$ (c.f.,~\cite{ horn1990matrix}) to obtain:
\begin{equation} \label{eq:sd}
K = MP\Lambda P^{-1},
\end{equation}
where $\Lambda$ is a diagonal matrix with entries to be the eigenvalues of the generalized eigenvalue problem
\begin{equation} \label{eq:eigKM}
K v = \Lambda Mv
\end{equation}
and $P$ is the matrix formed by the eigenvectors. We sort the eigenvalues in ascending order and are listed in $\Lambda$, and the $j$-th column of $P$ is associated with the eigenvalue $\lambda_j$ corresponds to the $j$-th diagonal entry of $\Lambda$.

Using~\eqref{eq:sd} and~\eqref{eq:G}, and considering $ I = P I P^{-1}$, we calculate
\begin{equation}
\begin{aligned}
G^{-1} & = (M + \eta K)^{-1} 
 = (M + \eta MP\Lambda P^{-1} )^{-1} \\
& = \Big( M P ( I + \eta \Lambda ) P^{-1} \Big)^{-1} 
 = P ( I + \eta \Lambda )^{-1} P^{-1} M^{-1}.
\end{aligned}
\end{equation}

Finally we obtain:
\begin{equation}
\begin{aligned}
G^{-1} K & = \Big( P ( I + \eta \Lambda )^{-1} P^{-1} M^{-1} \Big) \Big( MP\Lambda P^{-1} \Big) = P ( I + \eta \Lambda )^{-1} \Lambda P^{-1}, \\
G^{-1} M & = \Big( P ( I + \eta\Lambda )^{-1} P^{-1} M^{-1} \Big) M = P ( I + \eta \Lambda )^{-1} P^{-1}. \\
\end{aligned}
\end{equation}

If we define $E = ( I + \eta \Lambda )^{-1}$, we can rewrite the amplification matrix in~\eqref{eq:mp0} as:
\begin{equation}
\begin{aligned}
\Xi& =
\begin{bmatrix}
I - \tau^2 {\alpha_f}^2 J & I-\tau^2\alpha_f^3J& \frac{1}{2}I-\alpha_f^2R \\[0.2cm]
-\tau^2  \gamma J&I-\tau^2\gamma\alpha_fJ & I-\gamma R \\[0.2cm]
-\tau^2 J&-\tau^2 \alpha_fJ&I-R\\
\end{bmatrix}\\
 &=
\begin{bmatrix}
P & \bfs{0}&\bfs{0} \\
\bfs{0} & P &\bfs{0}\\
\bfs{0}&\bfs{0}&P\\
\end{bmatrix}
\begin{bmatrix}
I - \tau^2 \frac{{\alpha_f}^2}{\alpha_m} E\Lambda & I-\tau^2\frac{{\alpha_f}^2}{\alpha_m}E\Lambda& \frac{1}{2}+\frac{{\alpha_f}^2}{\alpha_m} \big(-E-\frac{\tau^2}{2}\alpha_fE\Lambda \big) \\[0.2cm]
-\tau^2 \frac{ \gamma }{\alpha_m} E\Lambda &I-\tau^2\gamma\frac{{\alpha_f}}{\alpha_m} E\Lambda &I+\frac{{\gamma}}{\alpha_m} \big(-E-\frac{\tau^2}{2}\alpha_fE\Lambda \big) \\[0.2cm]
-\frac{{\tau}^2}{\alpha_m} E\Lambda&-\tau^2\frac{{\alpha_f}}{\alpha_m}E\Lambda&I+\frac{1}{\alpha_m} \big(-E-\frac{\tau^2}{2}\alpha_fE\Lambda \big)\\
\end{bmatrix}\\
&\begin{bmatrix}
P^{-1} & \bfs{0}&\bfs{0} \\
\bfs{0} & P^{-1}&\bfs{0} \\
\bfs{0}&\bfs{0}&P^{-1}\\
\end{bmatrix}.\\
\end{aligned}
\end{equation}
Thus, we have
\begin{equation} \label{eq:mp1}
\begin{aligned}
&\begin{bmatrix}
U^n \\
\tau V^n\\
\tau^2 A^n\\
\end{bmatrix}
=
\begin{bmatrix}
P & \bfs{0}&\bfs{0} \\
\bfs{0} & P &\bfs{0}\\
\bfs{0}&\bfs{0}&P\\
\end{bmatrix}
\begin{bmatrix}
I - \tau^2 \frac{{\alpha_f}^2}{\alpha_m} E\Lambda & I-\tau^2\frac{{\alpha_f}^2}{\alpha_m}E\Lambda& \frac{1}{2}+\frac{{\alpha_f}^2}{\alpha_m} \big(-E-\frac{\tau^2}{2}\alpha_fE\Lambda \big) \\[0.2cm]
-\tau^2 \frac{ \gamma }{\alpha_m} E\Lambda &I-\tau^2\gamma\frac{{\alpha_f}}{\alpha_m} E\Lambda &I+\frac{{\gamma}}{\alpha_m} \big(-E-\frac{\tau^2}{2}\alpha_fE\Lambda \big) \\[0.2cm]
-\frac{{\tau}^2}{\alpha_m} E\Lambda&-\tau^2\frac{{\alpha_f}}{\alpha_m}E\Lambda&I+\frac{1}{\alpha_m} \big(-E-\frac{\tau^2}{2}\alpha_fE\Lambda \big)\\
\end{bmatrix}^n\\
&\begin{bmatrix}
P^{-1} & \bfs{0}&\bfs{0} \\
\bfs{0} & P^{-1}&\bfs{0} \\
\bfs{0}&\bfs{0}&P^{-1}\\
\end{bmatrix}
\begin{bmatrix}
U^0 \\
\tau V^0\\
\tau^2 A^0\\
\end{bmatrix}.
\end{aligned}
\end{equation}
Let us denote the matrix raised to power $n$ by $\tilde \Xi$. The method is unconditionally stable when the spectral radius of this matrix is bound by one. Herein, we omit the analysis for brevity and state that the method is unconditionally stable for specific values for $\alpha_m$ and $\alpha_f$ given in~\eqref{eq:ro}. We also refer the reader to~\cite{behnoudfar2018variationally} in which we detailed the analysis.

\subsection{Stability of the splitting schemes}\label{stb}

Here, we first analyze the scenario of only applying the approximating $\tilde{G}$ and then the proposed splitting approach. Similarly, we consider the spectral decomposition of each of the directional matrices $K_\xi$ with respect to its directional $M_\xi$ and obtain
\begin{equation} \label{eq:sd0}
K_\xi = M_\xi P_\xi D_\xi P_\xi^{-1},
\end{equation}
where $D_\xi$ is a diagonal matrix with entries being the eigenvalues of the generalized eigenvalue problem
\begin{equation} \label{eq:eigKM0}
K_\xi v_\xi = \lambda_\xi M_\xi v_\xi
\end{equation}
and $P_\xi$ is a matrix, with all the columns being the eigenvectors. Herein, $\xi = x,y,z$ specifies each of the coordinate directions. We state the analysis for 2D splitting and similarly to the case of generalized-$\alpha$ scheme; we calculate the required terms as follows (for more details see~\cite{ behnoudfar2018variationally})
\begin{equation}\label{eq:multi}
\begin{aligned}
\tilde G^{-1} &= P_x E_x P_x^{-1} M_x^{-1} \otimes P_y E_y P_y^{-1} M_y^{-1},\\
\tilde G^{-1} M &= \Big( P_x \otimes P_y \Big) \cdot \Big( E_x \otimes E_y \Big) \cdot \Big( P_x^{-1} \otimes P_y^{-1} \Big), \\
\tilde G^{-1} K & = \Big( P_x \otimes P_y \Big) \cdot \Big( E_x D_x \otimes E_y + E_x \otimes E_y D_y \Big) \cdot \Big( P_x^{-1} \otimes P_y^{-1} \Big). \\
\end{aligned}
\end{equation}
where:
\begin{equation}
E_\xi = ( I + \eta D_\xi )^{-1}, \qquad \xi = x, y.
\end{equation}

If we use the following identity:
\begin{equation}
\begin{aligned}
I & = P_x I_x P_x^{-1} \otimes P_y I_y P_y^{-1}  = \Big( P_x \otimes P_y \Big) \cdot \Big( I_x \otimes I_y \Big) \cdot \Big( P_x^{-1} \otimes P_y^{-1} \Big), \\
\end{aligned}
\end{equation}
then, the blocks of the amplification matrix are:
\begin{equation}
\begin{aligned}
\Xi_{11} & = \Big( P_x \otimes P_y \Big) \cdot \Big( I_x\otimes I_y - \tau^2 \frac{{\alpha_f}^2}{\alpha_m} \big(E_x D_x \otimes E_y + E_x \otimes E_y D_y\big) \Big) \cdot \Big( P_x^{-1} \otimes P_y^{-1} \Big), \\
\Xi_{12} &  = \Big( P_x \otimes P_y \Big) \cdot \Big(I_x\otimes I_y-\tau^2\frac{{\alpha_f}^3}{\alpha_m}\big(E_x D_x \otimes E_y + E_x \otimes E_y D_y\big) \Big) \cdot \Big( P_x^{-1} \otimes P_y^{-1} \Big), \\
\Xi_{13} &  = \Big( P_x \otimes P_y \Big) \cdot \Big( (\frac{1}{2}-\alpha_f^2)(I_x\otimes I_y)+\frac{{\alpha_f}^2}{\alpha_m} \big(\alpha_m(I_x\otimes I_y)-E\\&-\frac{\tau^2}{2}\alpha_f\big(E_x D_x \otimes E_y + E_x \otimes E_y D_y\big)\big) \Big) \cdot \Big( P_x^{-1} \otimes P_y^{-1} \Big), \\
\Xi_{21} &  = \Big( P_x \otimes P_y \Big) \cdot \Big( -\tau^2 \frac{\gamma}{\alpha_m}\big(E_x D_x \otimes E_y + E_x \otimes E_y D_y\big) \Big) \cdot \Big( P_x^{-1} \otimes P_y^{-1} \Big), \\
\Xi_{22} &  = \Big( P_x \otimes P_y \Big) \cdot \Big(  I_x\otimes I_y-\tau^2\gamma\frac{{\alpha_f}}{\alpha_m}\big(E_x D_x \otimes E_y + E_x \otimes E_y D_y\big) \Big) \cdot \Big( P_x^{-1} \otimes P_y^{-1} \Big), \\
\Xi_{23} &  = \Big( P_x \otimes P_y \Big) \cdot \Big(  \big(I_x\otimes I_y-\frac{{\gamma}}{\alpha_m}\big(E_x \otimes E_y\big)-\frac{\tau^2 \gamma}{2 \alpha_m}\alpha_f\big(E_x D_x \otimes E_y + E_x \otimes E_y D_y\big)\big) \Big) \cdot \Big( P_x^{-1} \otimes P_y^{-1} \Big), \\
\Xi_{31} &  = \Big( P_x \otimes P_y \Big) \cdot \Big( -\frac{\tau^2}{\alpha_m} \big(E_x D_x \otimes E_y + E_x \otimes E_y D_y\big) \Big) \cdot \Big( P_x^{-1} \otimes P_y^{-1} \Big), \\
\Xi_{32} &  = \Big( P_x \otimes P_y \Big) \cdot \Big( -\tau^2 \frac{\alpha_f}{\alpha_m}\big(E_x D_x \otimes E_y + E_x \otimes E_y D_y\big) \Big) \cdot \Big( P_x^{-1} \otimes P_y^{-1} \Big), \\
\Xi_{33} &  = \Big( P_x \otimes P_y \Big) \cdot \Big( I_x\otimes I_y-\frac{1}{\alpha_m}\big(E_x \otimes E_y\big)-\frac{\tau^2}{2\alpha_m}\alpha_f\big(E_x D_x \otimes E_y + E_x \otimes E_y D_y\big) \Big) \cdot \Big( P_x^{-1} \otimes P_y^{-1} \Big).\\
\end{aligned}
\end{equation}
By denoting $\zeta= E_x D_x \otimes E_y + E_x \otimes E_y D_y $ and $\tilde{E}=E_x \otimes E_y $, we write the matrix as:
\begin{equation}\label{eq:stability}
\begin{aligned}
\Xi & =
\begin{bmatrix}
P_x\otimes P_y & \bfs{0}&\bfs{0} \\
\bfs{0} & P_x\otimes P_y  &\bfs{0}\\
\bfs{0}&\bfs{0}&P_x\otimes P_y \\
\end{bmatrix}
\begin{bmatrix}
I - \tau^2 \frac{{\alpha_f}^2}{\alpha_m} \zeta & I-\tau^2\frac{{\alpha_f}^3}{\alpha_m}\zeta& \frac{1}{2}I+\frac{{\alpha_f}^2}{\alpha_m} \big(-\tilde{E}-\frac{\tau^2}{2}\alpha_f\zeta \big) \\[0.2cm]
-\tau^2  \frac{\gamma}{\alpha_m} \zeta &I-\tau^2\gamma\frac{{\alpha_f}}{\alpha_m}\zeta &   I-\frac{\gamma}{\alpha_m}\tilde{E}-\frac{\tau^2\gamma}{2\alpha_m}\alpha_f\zeta \\[0.2cm]
-\frac{\tau^2}{\alpha_m} \zeta&-\tau^2\frac{{\alpha_f}}{\alpha_m}\zeta&I-\frac{1}{\alpha_m}\tilde{E}-\frac{\tau^2}{2\alpha_m}\alpha_f\zeta\\
\end{bmatrix}^n\\
&\begin{bmatrix}
 P_x^{-1} \otimes P_y^{-1} & \bfs{0}&\bfs{0} \\
\bfs{0} & P_x^{-1} \otimes P_y^{-1}&\bfs{0} \\
\bfs{0}&\bfs{0}& P_x^{-1} \otimes P_y^{-1}\\
\end{bmatrix}.
\end{aligned}
\end{equation}

The stability of the method follows the same logic as analysis of the generalized-$\alpha$ method by calculating the spectral radius of:
\begin{equation}
\tilde{\Xi}  =
\begin{bmatrix}
I - \tau^2 \frac{{\alpha_f}^2}{\alpha_m} \zeta &I-\tau^2\frac{{\alpha_f}^3}{\alpha_m}\zeta& \frac{1}{2}I+\frac{{\alpha_f}^2}{\alpha_m} \big(-\tilde{E}-\frac{\tau^2}{2}\alpha_f\zeta \big) \\[0.2cm]
-\tau^2  \frac{\gamma}{\alpha_m} \zeta &I-\tau^2\gamma\frac{{\alpha_f}}{\alpha_m}\zeta &   I-\frac{\gamma}{\alpha_m}\tilde{E}-\frac{\tau^2\gamma}{2\alpha_m}\alpha_f\zeta\\[0.2cm]
-\frac{\tau^2}{\alpha_m} \zeta&-\tau^2\frac{{\alpha_f}}{\alpha_m}\zeta&I-\frac{1}{\alpha_m}\tilde{E}-\frac{\tau^2}{2\alpha_m}\alpha_f\zeta\\
\end{bmatrix}.\\
\end{equation}
First, by defining $\sigma= \tau^2 D_\xi$, we consider the two limiting cases for $\sigma$: $\sigma\rightarrow 0$ and $\sigma \rightarrow \infty$. In the limit $\sigma\to 0$, since $D_\xi$ is diagonal, $E_\xi \to I$ and consequently, we have $\tau^2 \zeta \rightarrow 0$ and $E \to I$ . Hence, $\tilde \Xi$ becomes upper triangular and the eigenvalues are obtained as:
\begin{equation}
\lambda_1=\lambda_2=1, \qquad \lambda_3=1-\frac{1}{\alpha_m}. 
\end{equation}
Hence, due to the equal multiplicity with the dimension of the stiffness matrix in 2D, $K$, the following condition is required:
\begin{equation} \label{eq:am}
\alpha_m \ge \frac{1}{2}.
\end{equation}

\begin{figure}[!ht]	
	\centering\includegraphics[width=6.2cm]{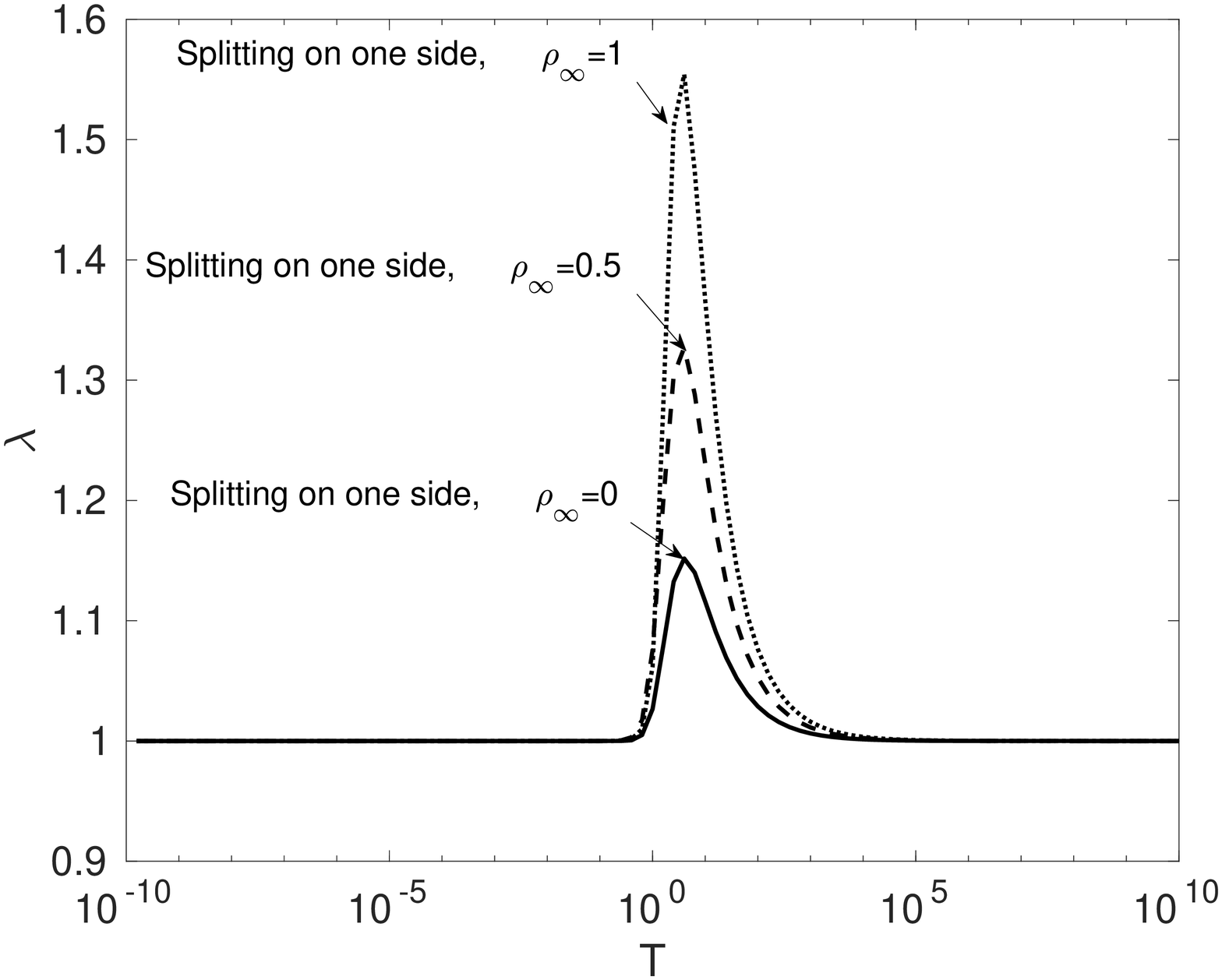} 
	\centering\includegraphics[width=6.5cm]{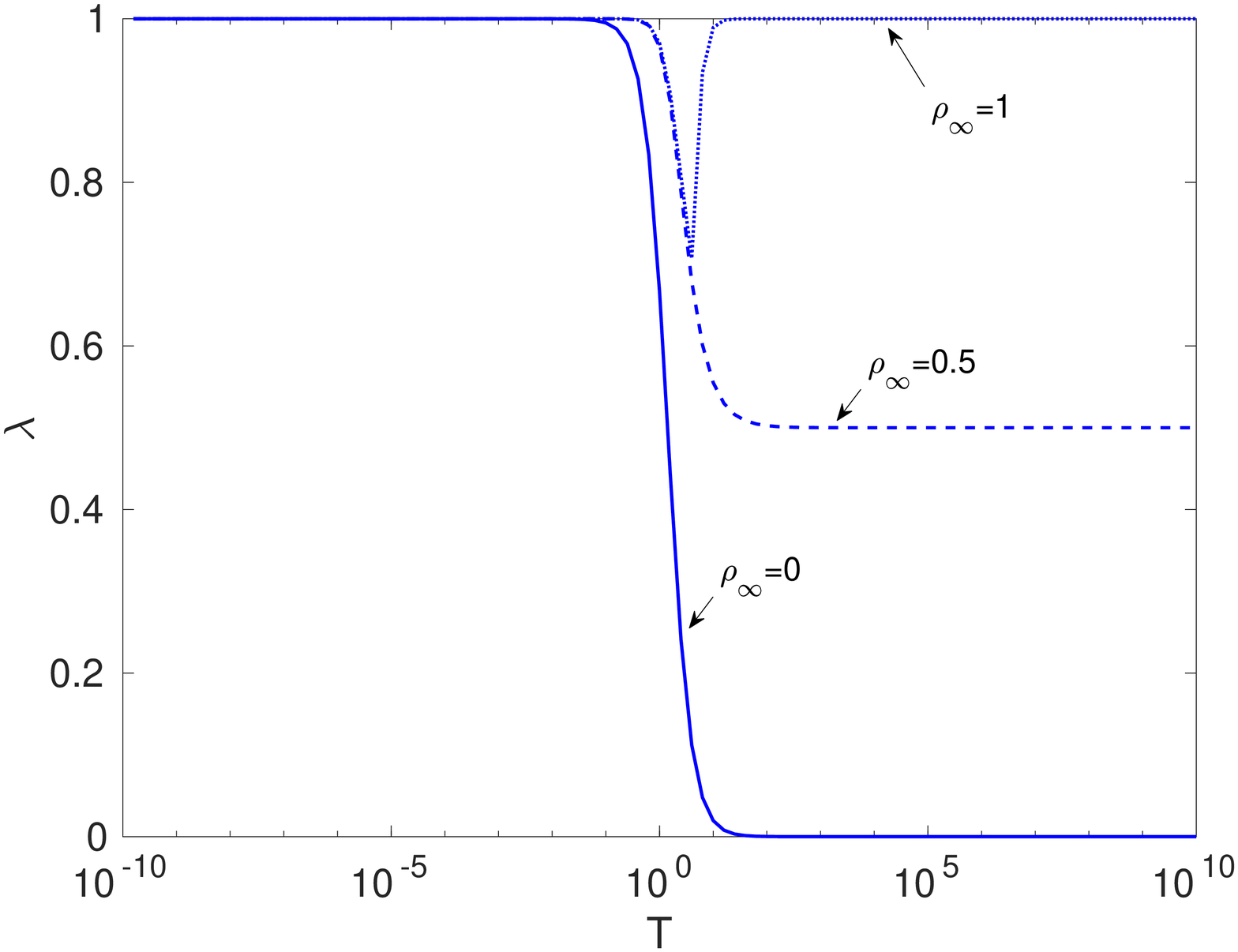} 	
	\caption{Largest eigenvalue of approximating $\tilde{G}$ (on the left) and the proposed splitting (on the right).}
	\label{fig:eig4}
	\end{figure}

In the case of infinite time step, the matrix $\tilde \Xi$ becomes:
\begin{equation}\label{eq:xi}
\begin{aligned}
\tilde \Xi& =
\begin{bmatrix}
I&I&\frac{1}{2}I\\
\bfs{0}&I&I\\
\bfs{0}&\bfs{0}&I\\
\end{bmatrix}.\\
\end{aligned}
\end{equation}
Here, we look into the proposed splitting method~\eqref{spli3} in more detail to discuss its need to consider the finite time-step sizes. We first show in figure~\ref{fig:eig4} that the largest eigenvalue for the system with $\tilde G$ only used on the left-hand side is not bounded by $1$ as $T$, scaled by $\tau^2 K$, grows. If we only split the left-hand side, the largest eigenvalue is not bounded for finite time steps, while it is bounded for the cases $\sigma \to 0$ or $\sigma \to \infty$. To address this unboundedness of the eigenvalue, we propose the splitting introduced in~\eqref{spli3}, which is unconditionally stable and provides dissipation control to the user even when $\sigma \to \infty$.
 By following a similar argument, we obtain the amplification matrix as:
\begin{equation}\label{eq:amp4}
\tilde{\Xi}  =
\begin{bmatrix}
I - \frac{1}{\alpha_f}\big(I-\tilde{E}\big)  &\tilde{E}& \frac{I}{2}-\beta \big(\frac{I}{2\beta}-(\frac{I}{2\beta}-1)\tilde{E}\big) \\[0.2cm]
 - \frac{\gamma}{\beta\alpha_f}\big(I-\tilde{E}\big) &I - \frac{\gamma}{\beta }\big(I-\tilde{E}\big) &   I-\gamma \big(\frac{I}{2\beta}-(\frac{I}{2\beta}-1)\tilde{E}\big) \\[0.2cm]
 - \frac{1}{\beta \alpha_f}\big(I-\tilde{E}\big) & - \frac{1}{\beta }\big(I-\tilde{E}\big) &I-\big(\frac{I}{2\beta}-(\frac{I}{2\beta}-1)\tilde{E}\big)\\
\end{bmatrix}.\\
\end{equation}
 For the case of $\sigma\to 0$, the eigenvalues of amplification matrix~\eqref{eq:amp4} are
\begin{equation}
\lambda_1=0, \qquad \qquad \lambda_2=\lambda_3=1.
\end{equation}
Likewise, for $\sigma \to \infty$, we have
\begin{equation}\label{eq:ampinf}
\tilde{\Xi}  =
\begin{bmatrix}
I - \frac{I}{\alpha_f}&\bfs{0} &\bfs{0}\\
-\frac{\gamma I}{\beta\alpha_f} &I - \frac{\gamma I}{\beta } &I-\frac{\gamma I}{2\beta}\\
-\frac{I}{\beta\alpha_f} & -\frac{I}{\beta} &I-\frac{I}{2\beta}\\
\end{bmatrix}.\\
\end{equation}
We show the bounded eigenvalues of~\eqref{eq:amp4} for the case of $T=1$, concluded from the figure~\ref{fig:eig4} as well as the eigenvalues of~\eqref{eq:ampinf} for various $\alpha_m$ and $\alpha_f$. 
\begin{figure}[!ht]	
	\centering\includegraphics[width=6.5cm]{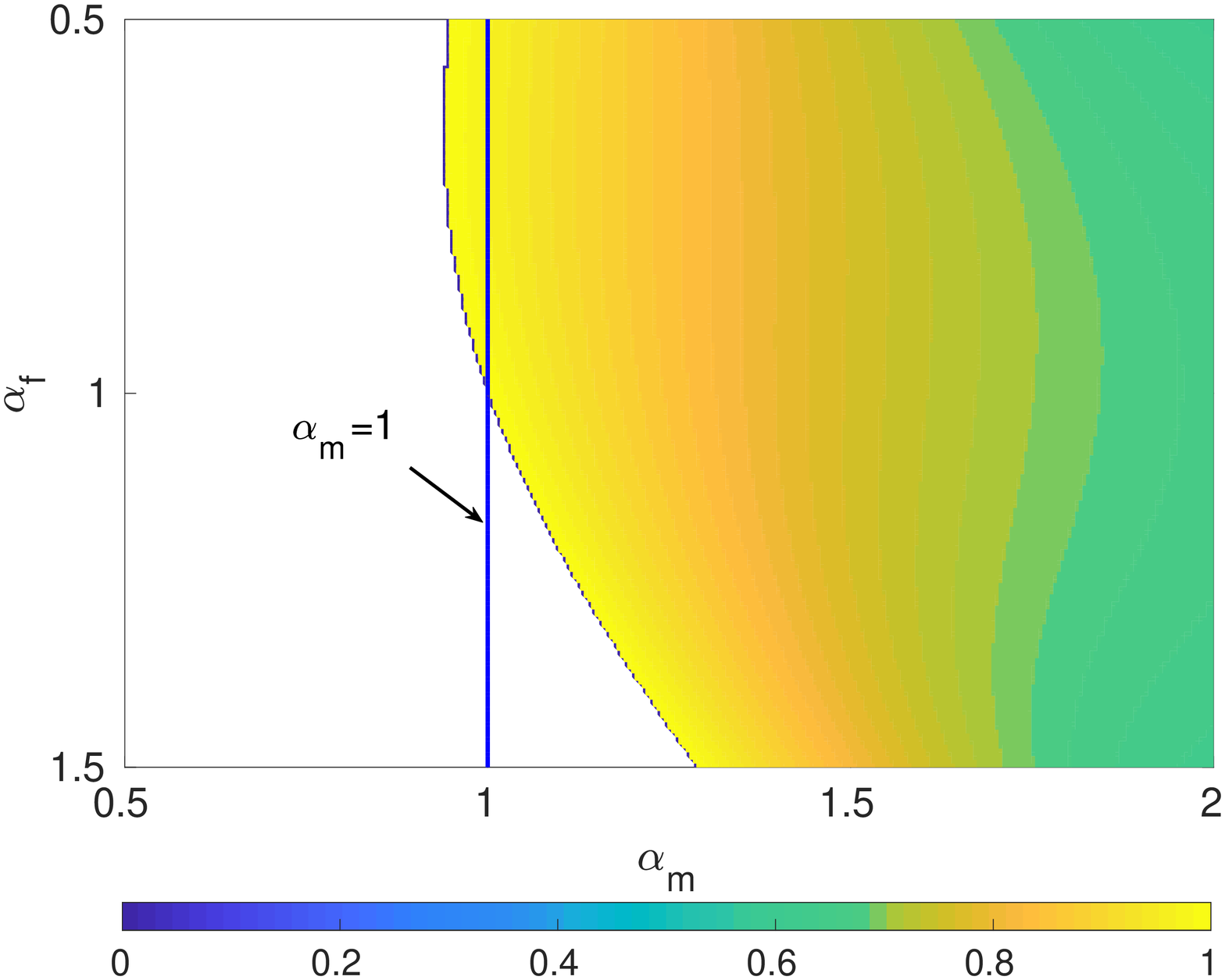} 
	\centering\includegraphics[width=6.5cm]{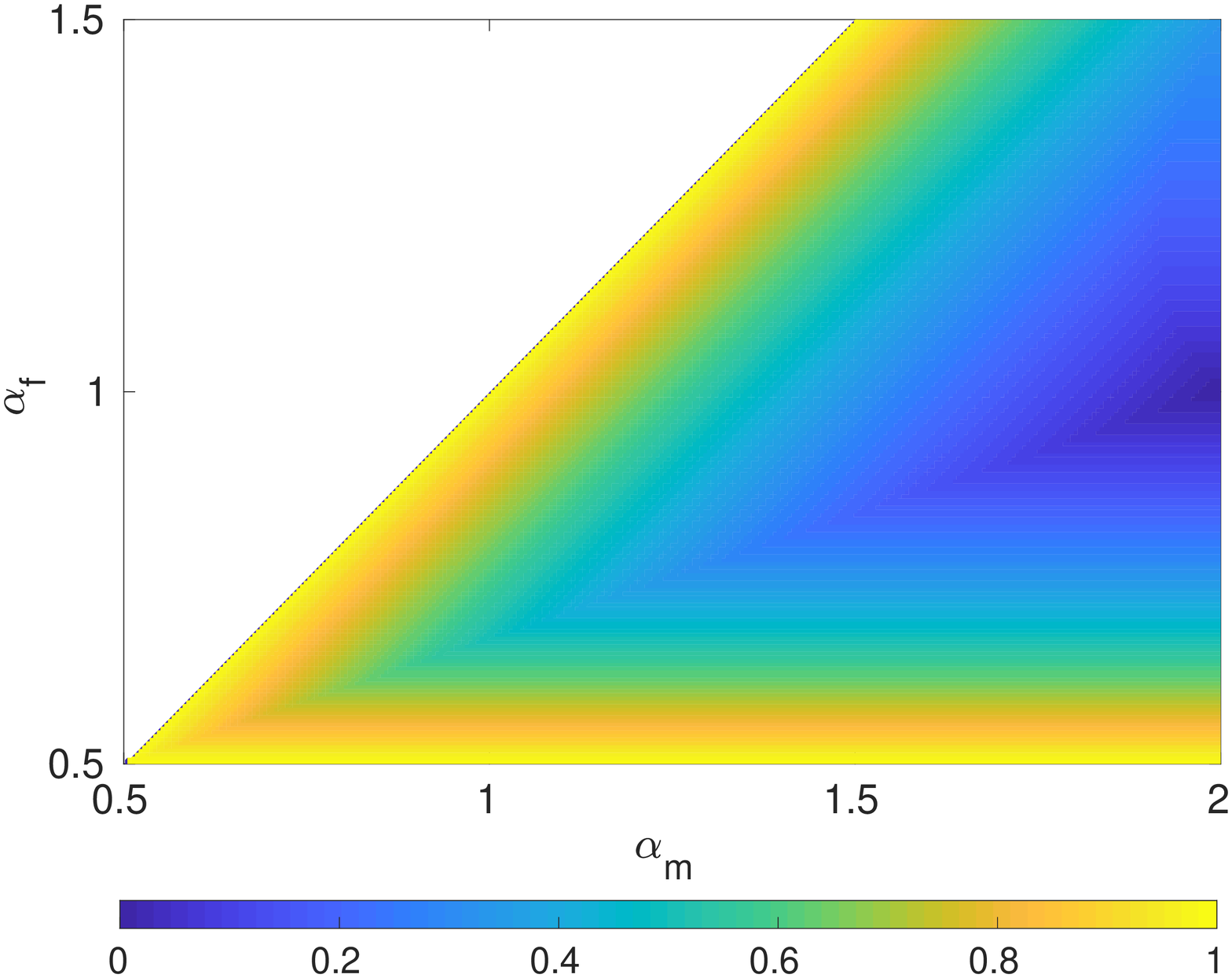} 	
	\caption{The region for bounded eigenvalues of the splitting method for $T=1$ (left) and $T \to \infty$ (right)}.
	\label{fig:eig}
\end{figure}
In order to control the high-frequency dissipation, we follow the idea proposed in~\cite{chung1993time}, by expressing the two parameters $\alpha_m$ and $\alpha_f$ in terms of the spectral radius $\rho_\infty$ of~\eqref{eq:ampinf}. By setting each of the eigenvalues of~\eqref{eq:ampinf} equal to $-\rho_\infty$, we state $	\alpha_f=\frac{1}{1+\rho_\infty}$ that is the similar to the generalized-$\alpha$ method. We also state that $\alpha_m=\frac{2-\rho_\infty}{1+\rho_\infty}$ for $\rho_\infty \leq 0.5$, and $\alpha_m=1$ for $0.5\le \rho_\infty \leq 1$. Hence, the method is unconditionally stable as well as A-stable by setting $\rho_\infty=0$.

\color{black}
\begin{remark}
	The conditions on $\alpha_m, \alpha_f$ for unconditional stability for the 3D splitting are the same as for the 2D splitting. The analysis is more involved, but it follows the same logic. We do not include it for brevity.
\end{remark}
Finally, we present the corresponding error estimates. From the above analysis, the splitting schemes have second-order accuracy in time.  Stability and consistency imply the method's convergence.  Thus, following the classical estimation (see, for example,~\cite{ thomee1984galerkin}), for regular solutions, we have the error estimates
\begin{equation} \label{eq:errl2h1}
\begin{aligned}
\| u_h^n - u(t_n) \|_{0, \Omega} & \le C(u) (h^{p+1} + \tau^2), \\
| u_h^n - u(t_n) |_{1, \Omega} & \le C(u) (h^p + \tau^2), \\
\end{aligned}
\end{equation}
where $u_h^n$ is the approximate solution at time $t_n$ and $C(u)$ is a positive constant independent of the mesh size $h$ and the time-step size $\tau$. 

\section{Numerical examples} \label{sec:5}

We now use the splitting methods proposed on several numerical examples to validate the efficiency and accuracy of the proposed methods. In all cases, we obtain optimal convergence rates for the spatial and temporal resolutions of the mesh. Additionally, we show a linear computational cost of the splitting schemes with respect to the total number of degrees of freedom in the system. Here, we consider the partial differential model problem~\eqref{eq:pde} with the corresponding forcing function, boundary, and initial conditions derived from an exact solution
\begin{equation}\label{pro}
u = 
\begin{cases}
u(x,y,t)=\sin(\pi x)\sin(\pi y) \Big (\sin (\sqrt{2} \pi t)+\cos (\sqrt{2} \pi t) \Big), \quad \text{in} \quad 2D, \\
u(x,y,z,t)=\sin(\pi x)\sin(\pi y) \sin(\pi z)  \Big (\sin (\sqrt{3} \pi t)+\cos (\sqrt{3} \pi t) \Big), \quad \text{in} \quad 3D. \\
\end{cases}
\end{equation}

We first validate our linear computation cost estimates. The inversion of the matrix is the main cost of solving~\eqref{eq:wave}. Here, the required number of operations for a 1D case is $\mathcal{O} (m)$ as a function of the degrees of freedom $m$ when using Gaussian elimination to invert the matrix, $M+\eta K$. By adopting the proposed splitting techniques, inversion of $\tilde{G}$ requires $\mathcal{O}(m_xm_y)$ and $\mathcal{O}(m_xm_ym_z)$ for the 2D and 3D cases, respectively. 

\begin{figure}[!ht]	
	\centering\includegraphics[width=6.5cm]{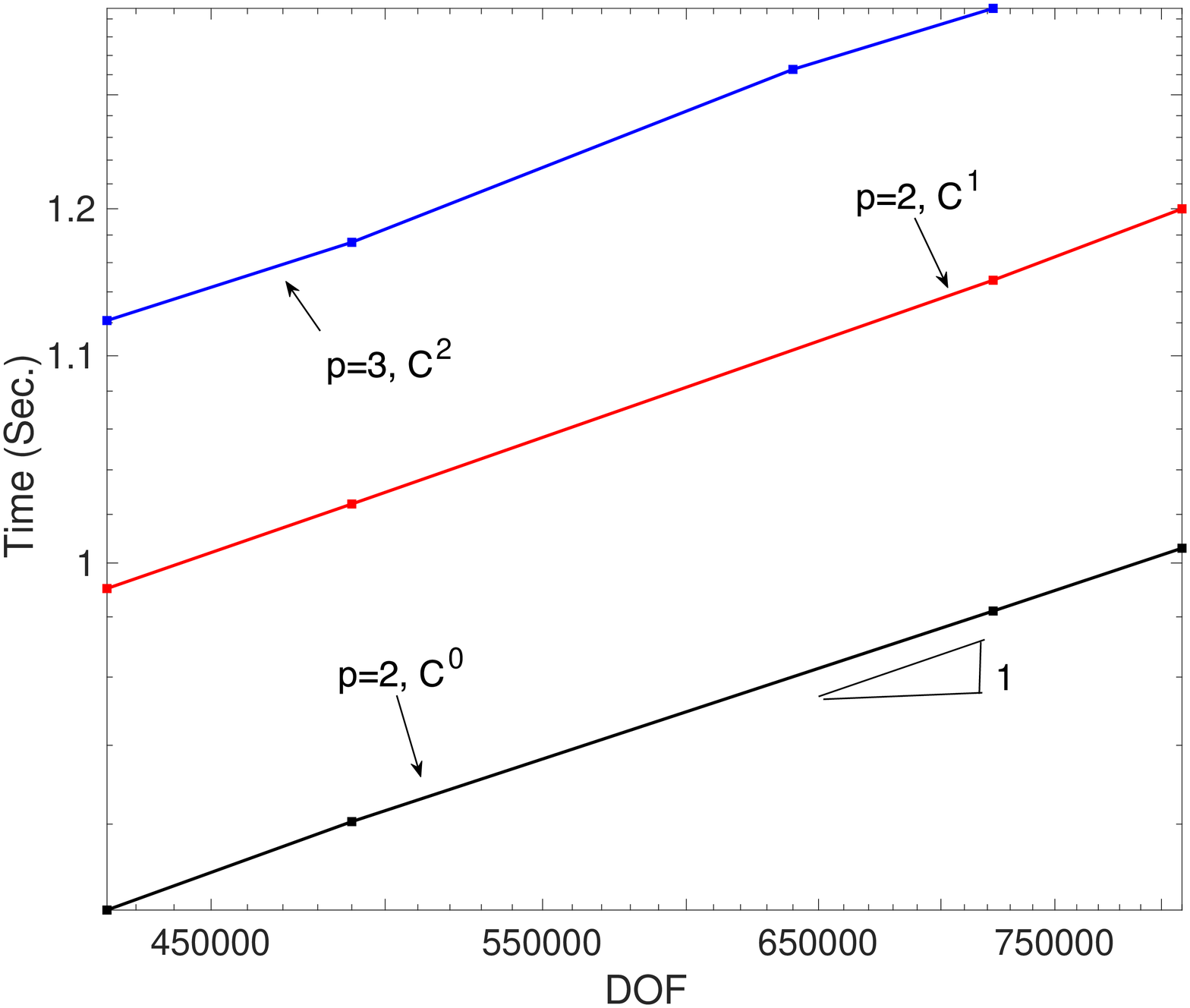} 
	\centering\includegraphics[width=6.5cm]{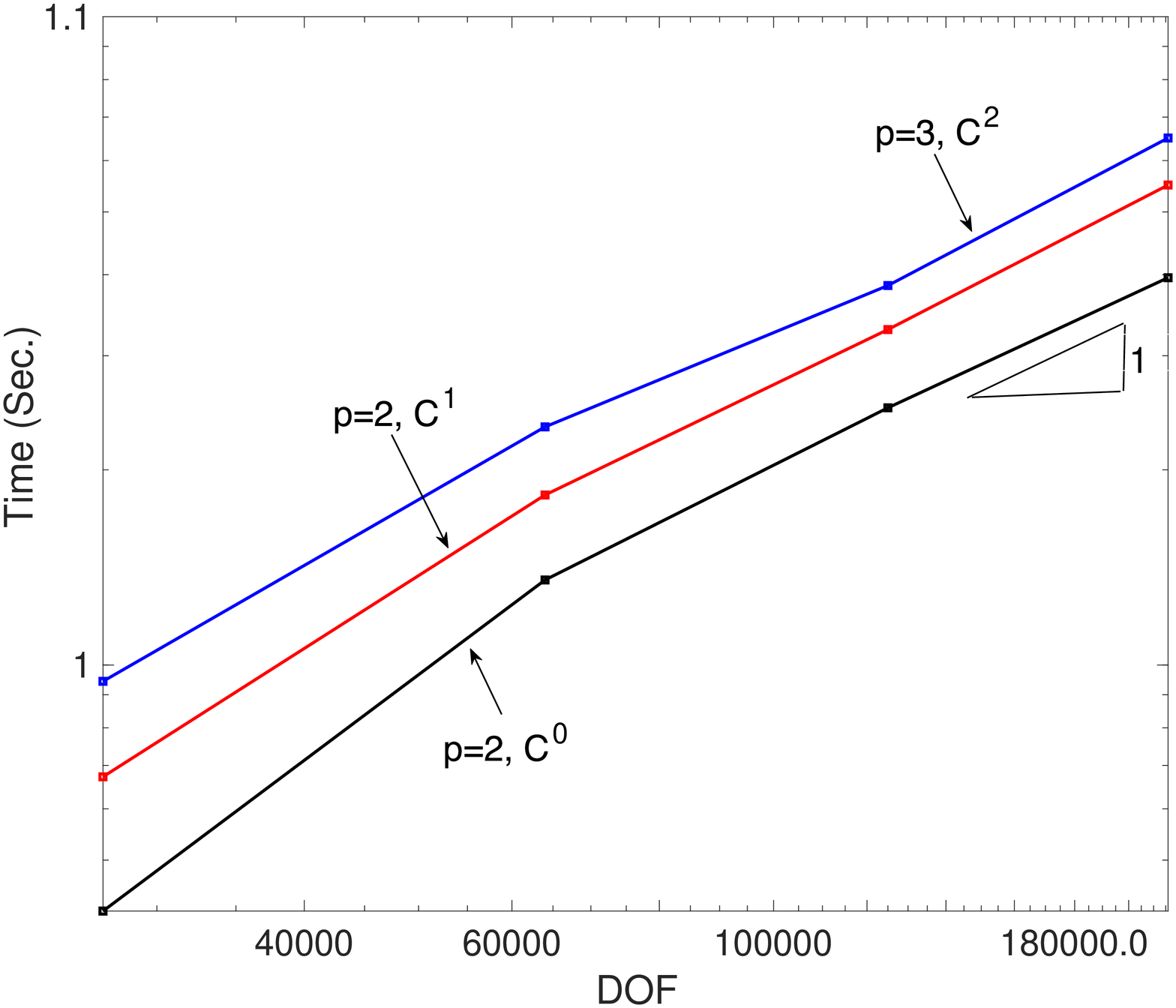} 	
\caption{The computational cost (linear) of the proposed splitting scheme when using $C^0$ and $C^1$ quadratic elements and $C^2$ cubic isogeometric elements with $\rho_\infty=0.5$ in 2D (left) and in 3D (right).}
	\label{fig:cost}
\end{figure}

To show the computational cost for the wave propagation problem~\eqref{eq:pde} in both 2D and 3D, we refer to figure~\ref{fig:cost}. This figure shows that the required cost for solving the multi-dimensional matrix problems of the proposed splitting scheme grows linearly to the total number of degrees of freedom in the system. Herein, as an example, we use a direct solver that is Gaussian elimination and three settings of the $C^0$ and $C^1$ quadratic elements as well as $C^2$ cubic isogeometric elements for the spatial discretization. The figure shows linear cost when applying the method for 2D and 3D cases. This proportionality validates the efficiency of the splitting scheme when solving the resulting matrix problems. Additionally, this approximation allows the use of direct solvers for problems of arbitrary dimension.

\begin{figure}[!ht]
	{\centering\includegraphics[width=6.5cm]{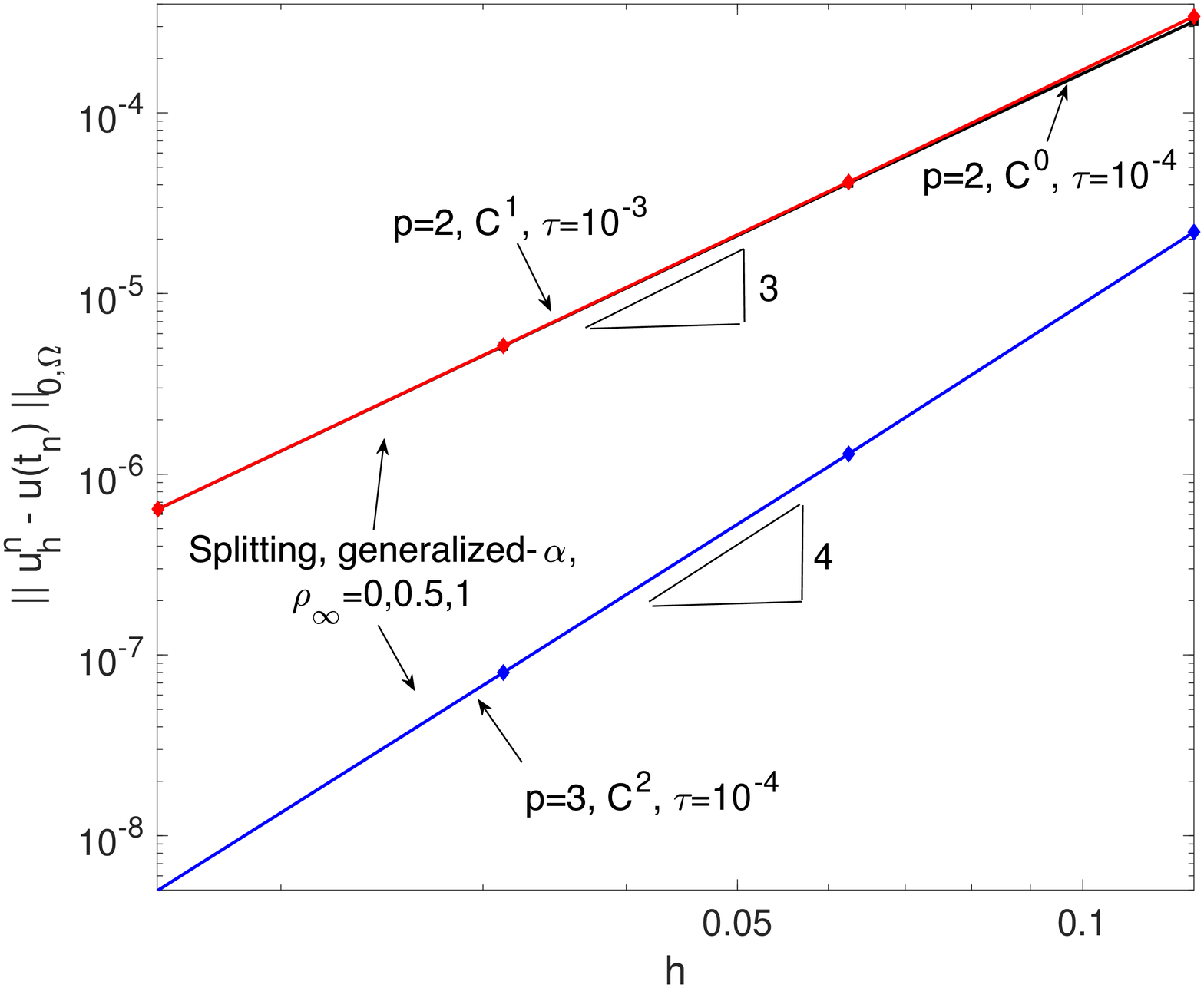}} 
		{\centering\includegraphics[width=6.5cm]{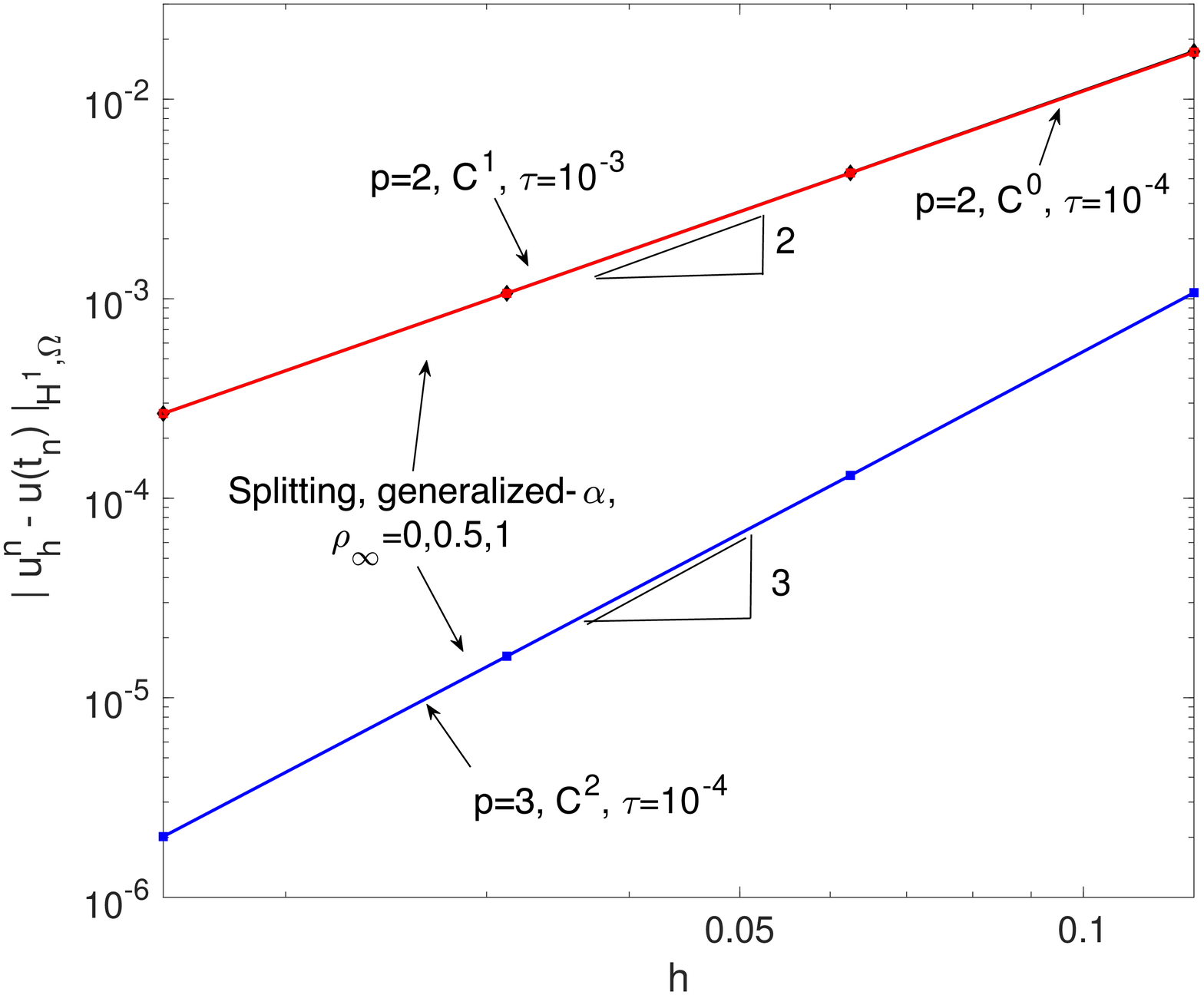}}
	\caption{Space-convergence of the splitting technique as well as generalized-$\alpha$ method for $\rho_\infty=0,0.5,1 $ in $L^2 $ norm and $H^1 $ semi-norm when using $p=2, C^1$, $p=3, C^2$ and classical finite elements $p=2, C^0$ for space discretization.}
	\label{fig:IGA}
\end{figure}

\begin{figure}[!ht]	
	{\centering\includegraphics[width=6.5cm]{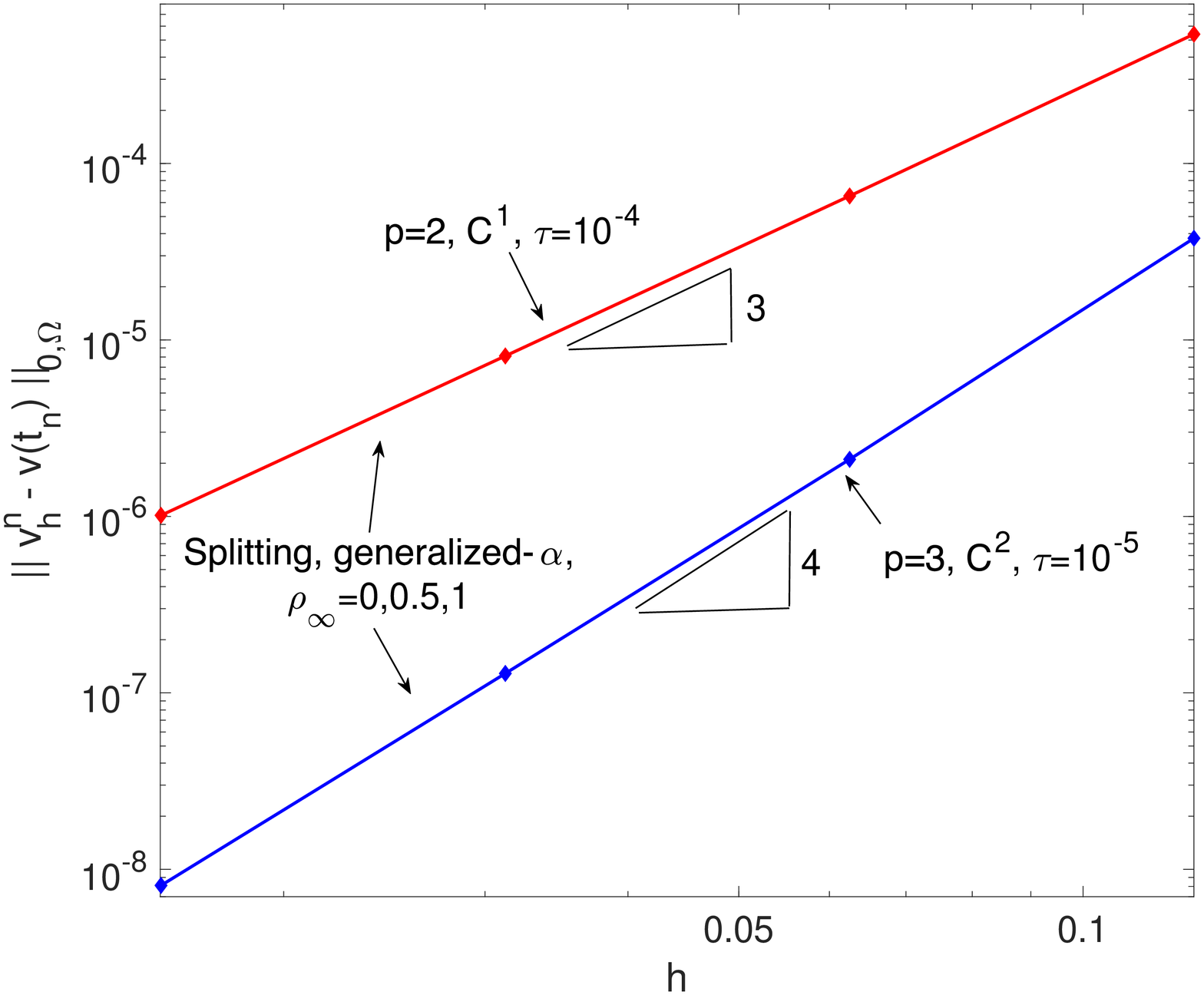}} 
{\centering\includegraphics[width=6.5cm]{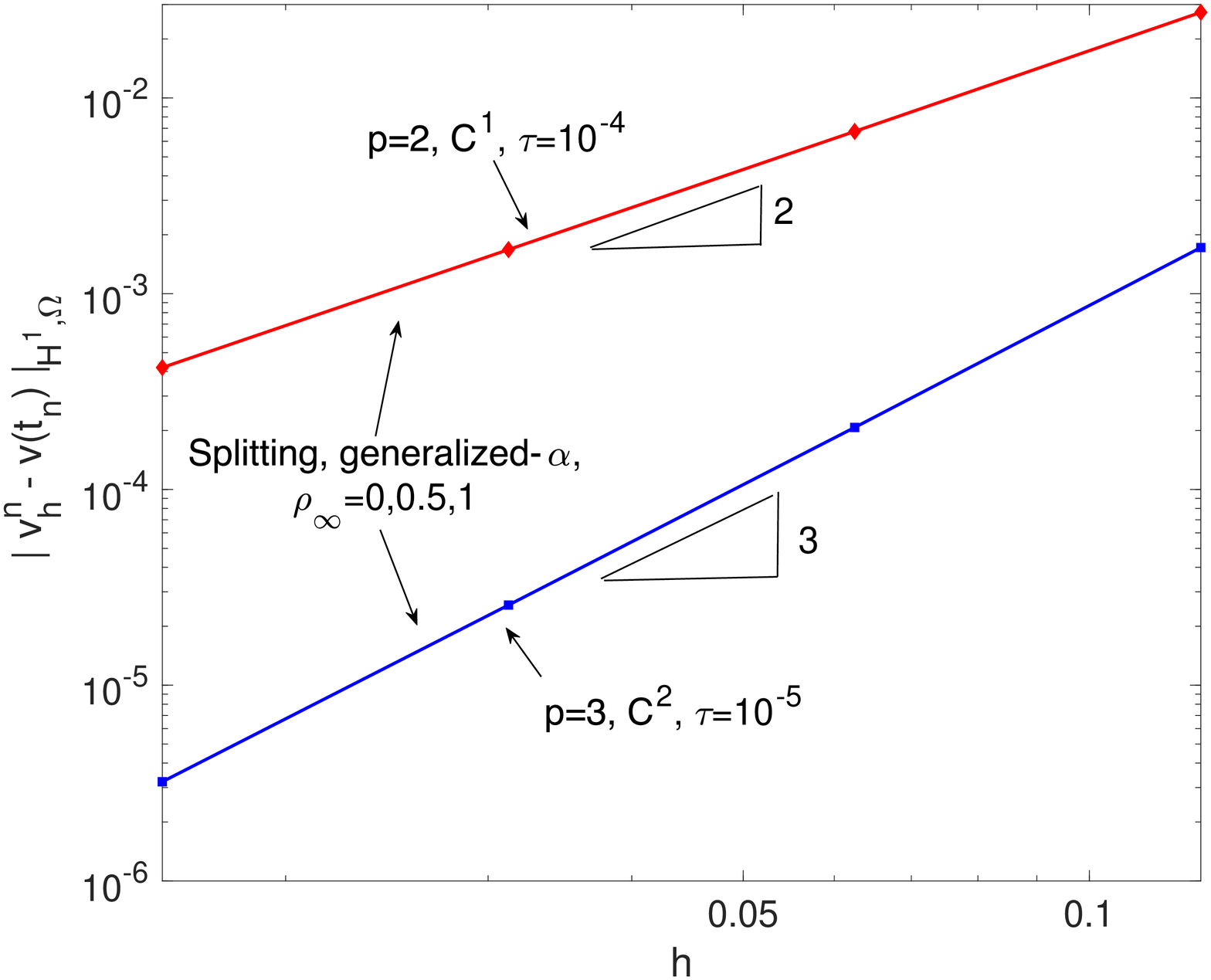}}
\caption{Space-convergence of the velocity obtained using the splitting technique and generalized-$\alpha$ method for $\rho_\infty=0,0.5,1 $ in $L^2 $ norm and $H^1 $ semi-norm when using $p=2, C^1$, $p=3, C^2$ for space discretization. }	\label{fig:FEM}
\end{figure}

In figure~\ref{fig:IGA}, we show the $L^2$ norm and $H^1$ semi-norm to study the space convergence of the methods for displacement for the final time of $0.1$ when choosing the fixed time steps $\tau=10^{-3}$ and $\tau=10^{-4} $ for $C^1$ quadratic and $C^2$ cubic elements, respectively. Additionally, we show that the method can also be used in conjunction with the classical finite element and delivers the optimal convergence rates, $p+1$ in $L^2$ norm, and $p$ in $H^1$ norm. Here, the splitting technique delivers the same error as the direct solution of the generalized-$\alpha$ method for $\rho_\infty=0,0.5,1$. Figure~\ref{fig:FEM} also shows the space convergence of the velocity for $C^1$ quadratic and $C^2$ cubic isogeometric elements.  

\begin{figure}[!ht]	
    \subfigure[Displacement: $p=2, C^0$]{\centering\includegraphics[width=6.5cm]{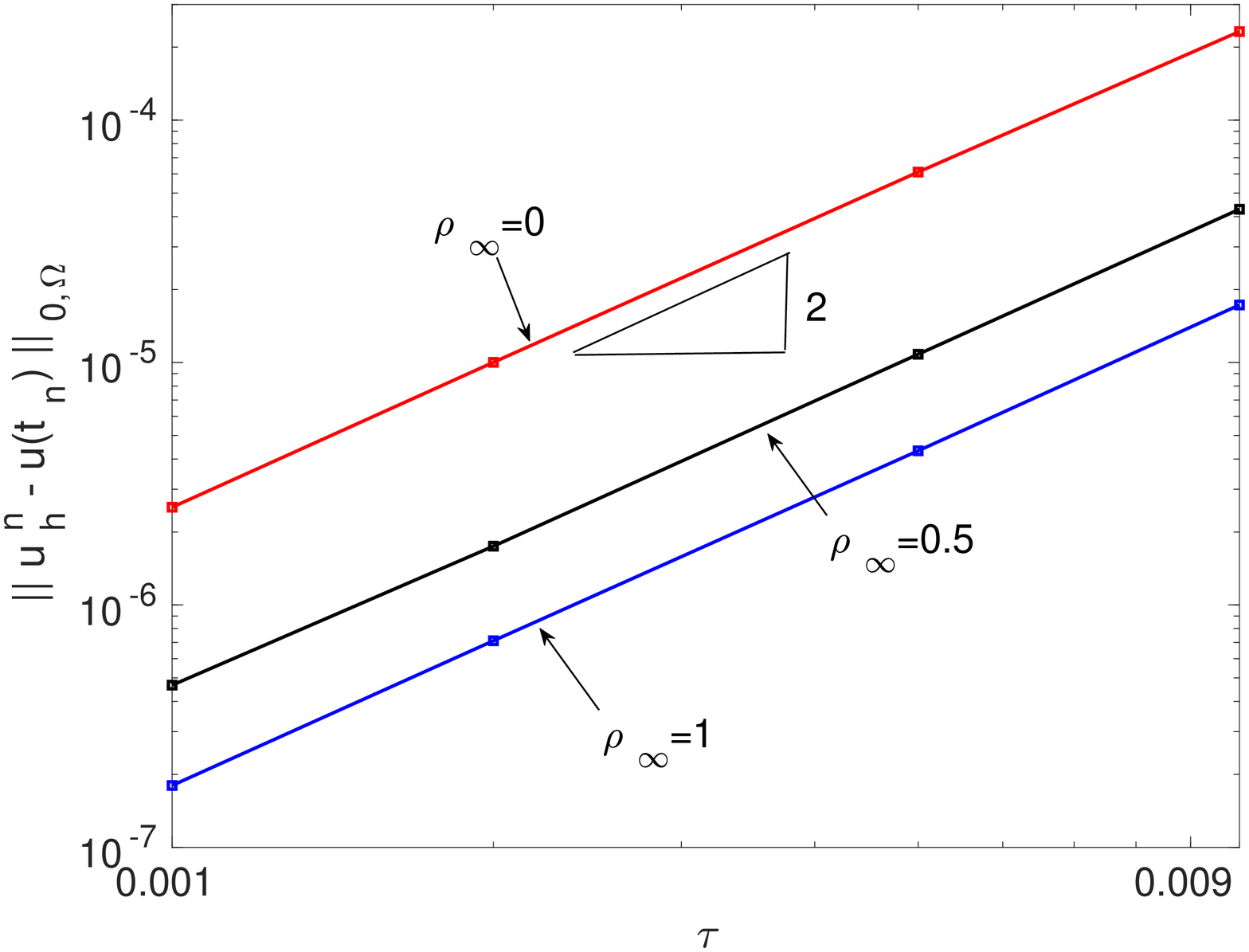}}
    	\subfigure[Velocity: $p=2, C^0$]{\centering\includegraphics[width=6.5cm]{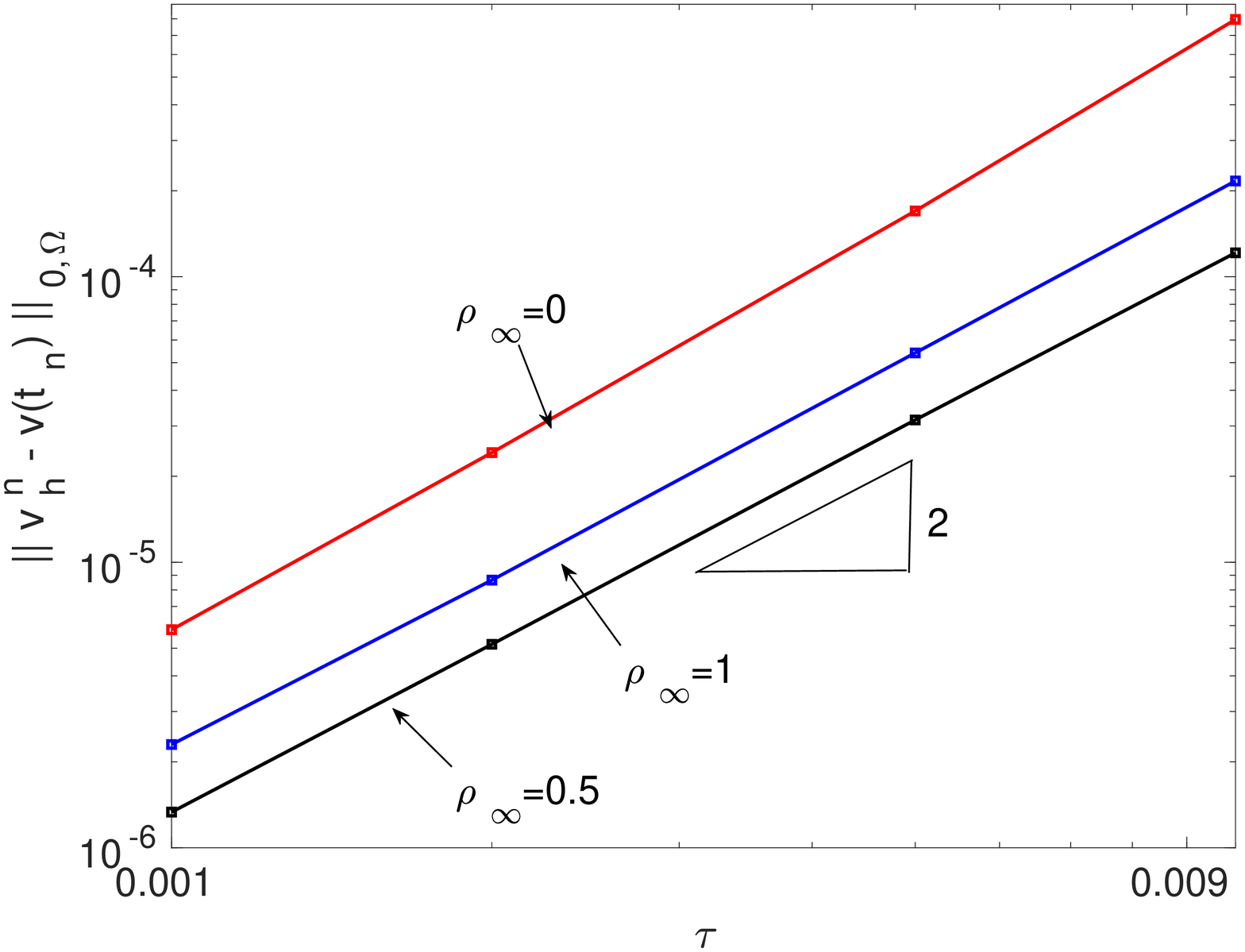}}
	\subfigure[Displacement: $p=2, C^1$]{\centering\includegraphics[width=6.5cm]{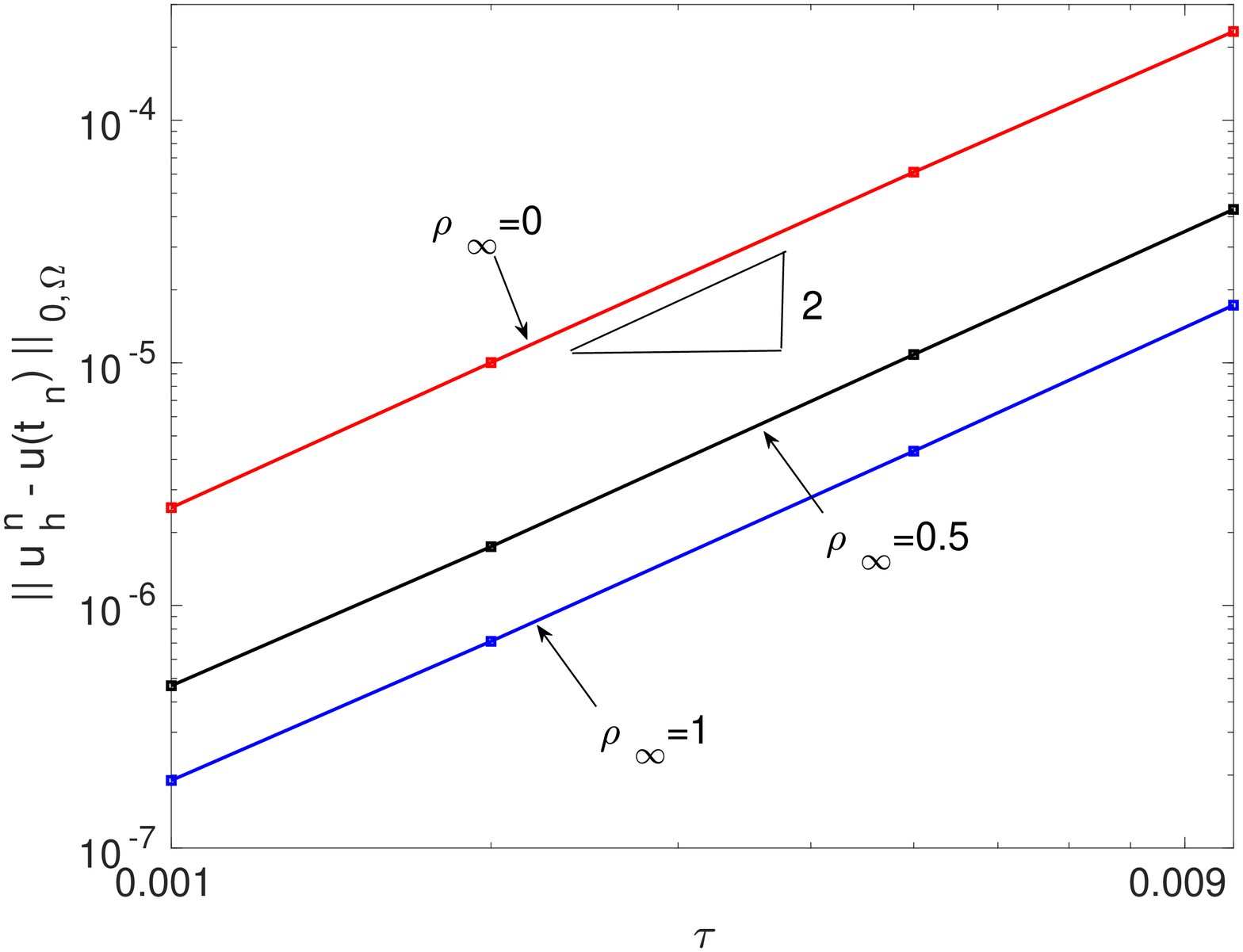}}
	\subfigure[Velocity: $p=2, C^1$]{\centering\includegraphics[width=6.5cm]{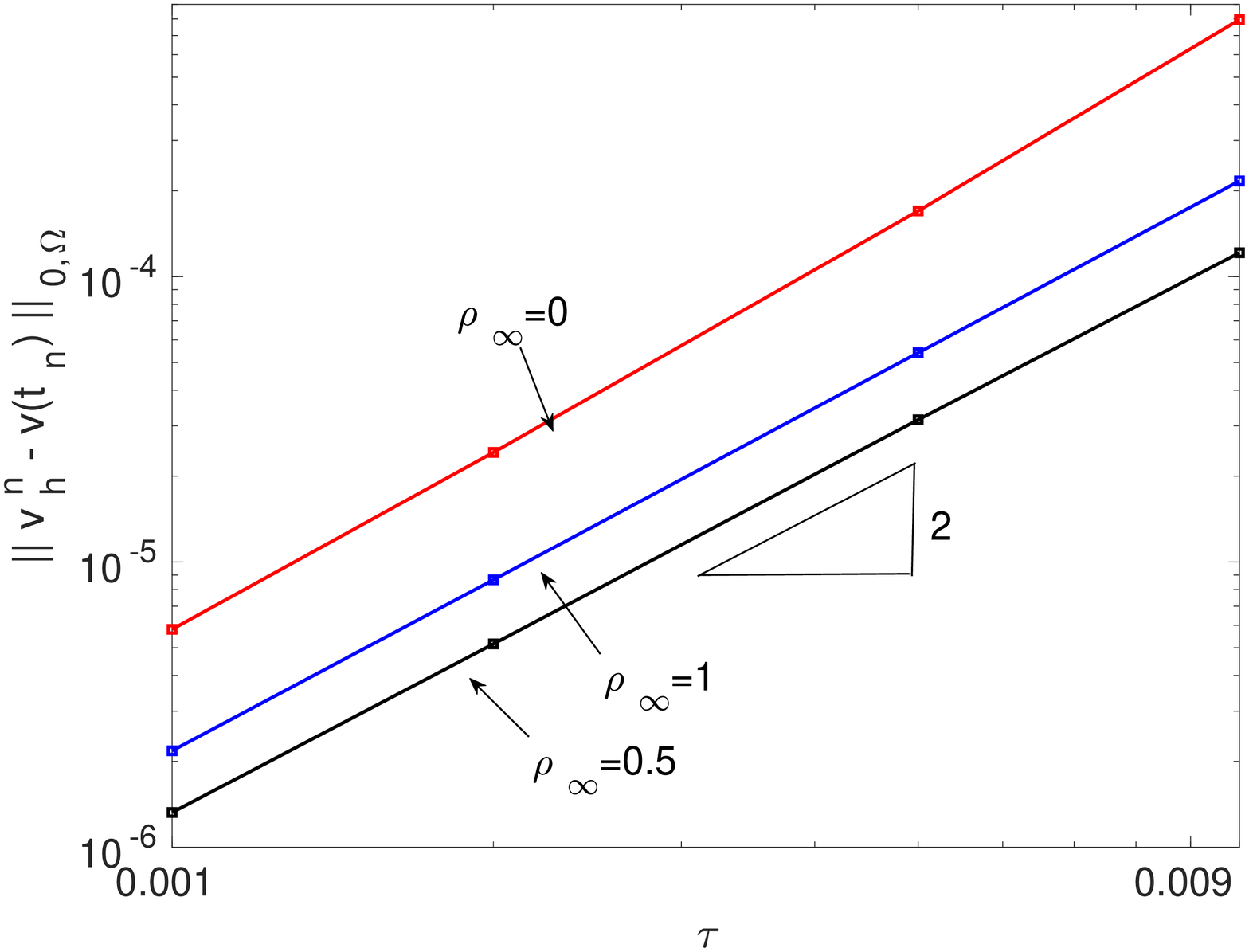}}
		\subfigure[Displacement: $p=3, C^2$]{\centering\includegraphics[width=6.5cm]{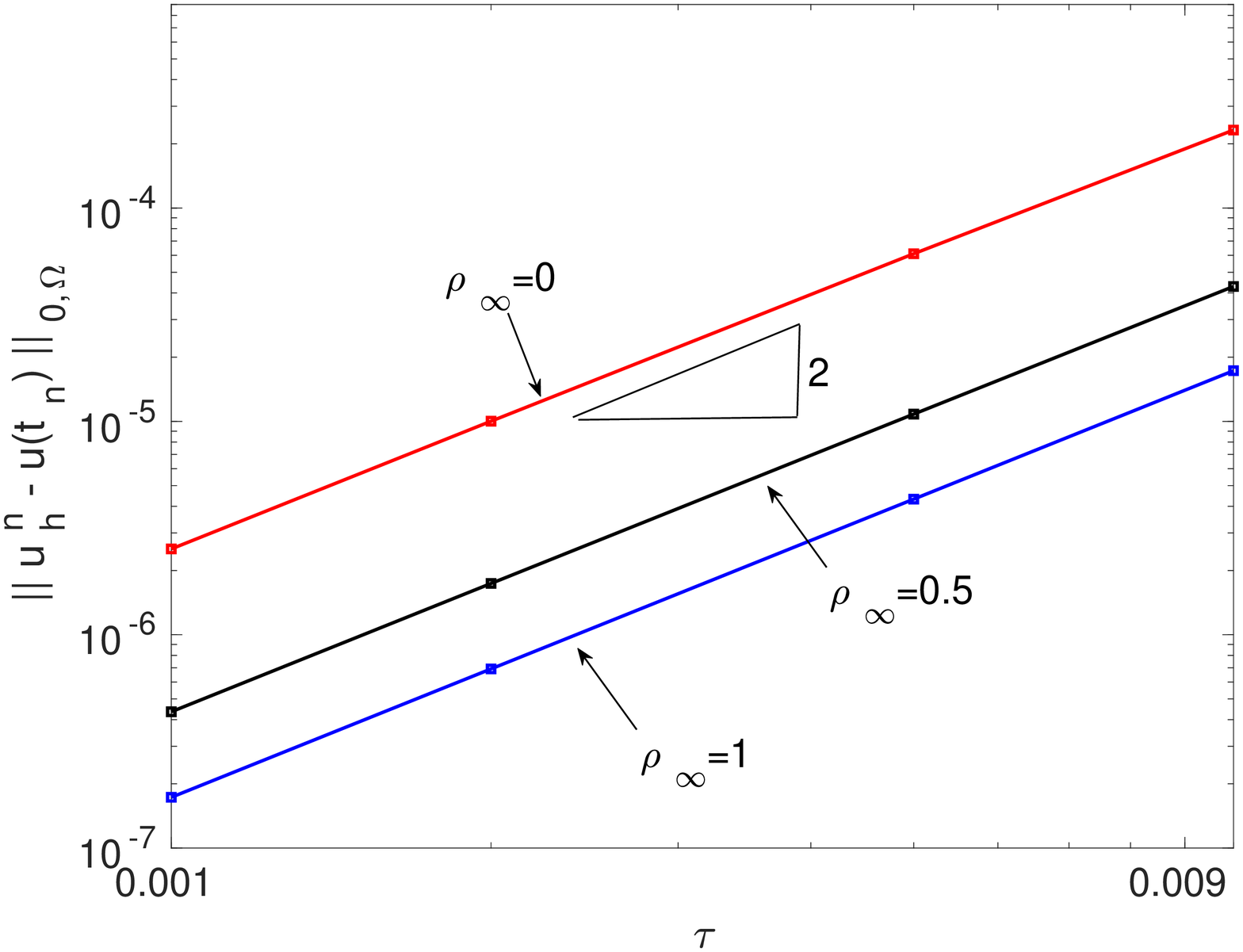}}
		\hspace{0.2 cm}
	\subfigure[Velocity: $p=3, C^2$]{\centering\includegraphics[width=6.5cm]{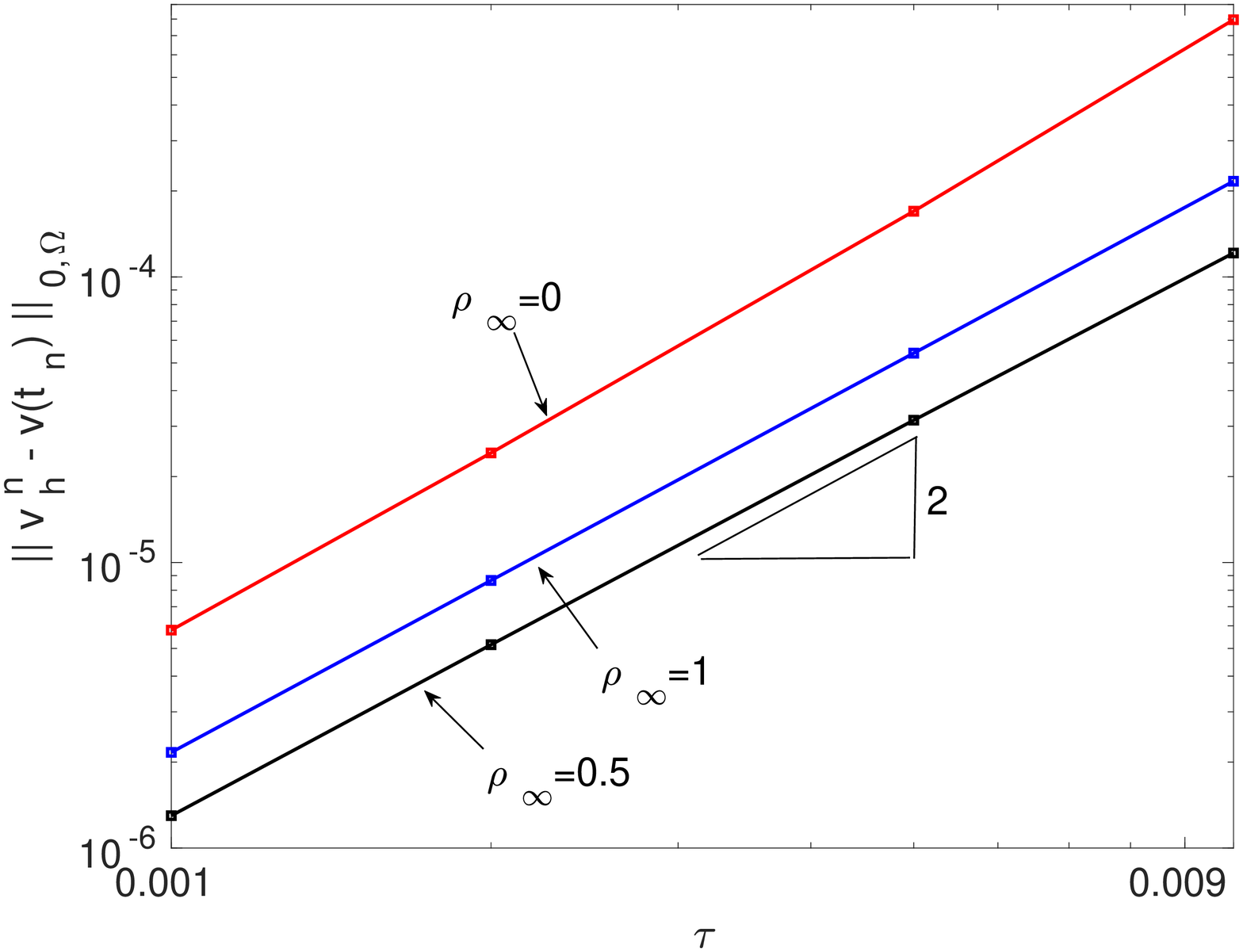}}
	\caption{Time-convergence of the solution obtained using the splitting method, in $L^2 $ norm using classical FEM: $p=2, C^0, n=100$ and IGA: $p=2, C^1, n=100$, and $p=3, C^2, n=64$ for space discretization.}
	\label{fig:time}
\end{figure}

Lastly, we provide numerical evidence of the second-order convergence
in time for our method. Figure~\ref{fig:time} illustrates $L^2$ norm for the case in which the number of elements $N=100\times100$ for FEM and IGA with quadratic elements, and $N=64\times64$ for IGA cubic elements. The final time $T$ is set to be $0.1$. 


\section{Concluding remarks}\label{sec:6}

We introduce a splitting technique for hyperbolic equations, which model wave propagation and structural dynamics. The method modifies the generalized-$\alpha$ method for temporal discretization.  The proposed splitting method is unconditionally stable and provides second-order accuracy in time as well as optimal convergence rates in space.  Another essential feature of the schemes is that the cost of solving the resulting algebraic system is proportional to the total number of degrees of freedom. Hence, the cost of the overall computation diminishes significantly for multidimensional problems.

These splitting schemes will be extended to higher-order generalized-$\alpha$ methods that we propose in~\cite{behnoudfar2019higher} and~\cite{deng2019high} as well as phase-field problems and to other time integrators. These advances will be reported shortly.

\section*{Acknowledgement}

This publication was made possible in part by the CSIRO Professorial Chair in Computational Geoscience at Curtin University and the Deep Earth Imaging Enterprise Future Science Platforms of the Commonwealth Scientific Industrial Research Organization, CSIRO, of Australia. Additional support was provided by the European Union’s Horizon 2020 Research and Innovation Program of the Marie Sklodowska-Curie grant agreement No. 777778, the Institute for Geoscience Research (TIGeR), and the Curtin Institute for Computation. The first and second authors also acknowledge the contribution of an Australian Government Research Training Program Scholarship in supporting this research.


\bibliographystyle{elsarticle-harv}\biboptions{square,sort,comma,numbers}
\bibliography{ref}

\end{document}